\documentclass[11pt]{article}
\usepackage{amssymb,amsfonts}
\usepackage[english]{babel}
\usepackage[pdftex]{graphicx}

\begin{document}

\begin{center}
{\large\bf An analytic proof of the Malgrange--Sibuya theorem on \\ the convergence of formal solutions of an ODE}
\end{center}

\begin{center}
R.\,R.\,Gontsov, I.\,V.\,Goryuchkina 
\end{center}
\bigskip

\begin{abstract} We propose an analytic proof of the Malgrange--Sibuya theorem concerning a sufficient condition of the convergence
of a formal power series satisfying an ordinary differential equation. The proof is based on the majorant method and allows to estimate
the radius of convergence of such a series.
\end{abstract}
\bigskip

{\bf\S1. Introduction}
\\

In the paper we study some properties of formal power series satisfying an ordinary differential equation
\begin{eqnarray}\label{ADE}
F(z,u,u',\ldots,u^{(n)})=0
\end{eqnarray}
of order $n$, where $F\not\equiv0$ is a holomorphic (in some domain) function of $n+2$ variables.

According to Maillet's theorem \cite{Mai}, if a formal series
$\hat\varphi=\sum_{j=0}^{\infty}c_jz^j\in{\mathbb C}[[z]]$
satisfies the equation (\ref{ADE}), where $F$ is a polynomial, then there is a real number
$s\geqslant0$ such that the power series $\sum_{j=0}^{\infty}(c_j/(j!)^s)z^j$
converges in some neighbourhood of zero. In this case one says that the formal series
$\hat\varphi$ has the {\it Gevrey type of order $s$}. Furthermore if there is no a real number
$s'<s$ such that the series has the Gevrey type of order $s'$, then the order $s$ is called {\it precise}.

How small the number $s$ from Maillet's theorem can be? And when one can guarantee the convergence
of the formal solution $\hat\varphi$ in a neighbourhood of zero, that is, when $s=0$? Answers on
these questions have been obtained in the papers of J.-P.\,Ramis, B.\,Malgrange, Y.\,Sibuya.

Let us first consider the case of a {\it linear} differential equation
\begin{eqnarray}\label{LDE}
b_n(z)u^{(n)}+b_{n-1}(z)u^{(n-1)}+\ldots+b_0(z)u=0,
\end{eqnarray}
whose coefficients $b_i(z)$ are holomorphic functions in a neighbourhood of zero. The point $z=0$ is said
to be a {\it Fuchsian} ({\it regular}) singular point of this equation, if its coefficients satisfy the conditions
$$
{\rm ord}_0\,b_{n-i}(z)+i\geqslant{\rm ord}_0\,b_n(z), \qquad
i=1,\ldots,n.
$$
The equation (\ref{LDE}) can be written in an equivalent form
\begin{eqnarray}\label{LDEdelta}
a_n(z)\delta^nu+a_{n-1}(z)\delta^{n-1}u+\ldots+a_0(z)u=0,
\end{eqnarray}
where $\displaystyle\delta=z\frac d{dz}$ and the coefficients $a_i(z)$ are holomorphic functions in a
neighbourhood of zero. With respect to such a form, the point $z=0$ is a Fuchsian singular point if the
inequalities
$$
{\rm ord}_0\,a_i(z)\geqslant{\rm ord}_0\,a_n(z), \qquad
i=0,\ldots,n-1,
$$
hold.

In the case when $z=0$ is a Fuchsian singular point of the equation (\ref{LDE}), any formal power
series that satisfies this equation converges in a neighbourhood of zero. Indeed, using a standard
substitution
$$
y^1=u, \qquad y^2=u', \qquad\ldots,\qquad y^n=u^{(n-1)},
$$
one can pass to a linear differential system with a {\it regular} singular point $z=0$, which is written
in a matrix form:
$$
\frac{dy}{dz}=B(z)\,y, \qquad y=(y^1,\ldots,y^n)^{\top}.
$$
This system, by means of a gauge transformation $\tilde y=\Gamma(z)\,y$ with a meromorphic matrix $\Gamma(z)$ at
zero, is equivalent to a (Fuchsian) system
\begin{equation}\label{System}
\frac{d\tilde y}{dz}=\frac Az\,\tilde y,
\end{equation}
where $A$ is a constant $(n\times n)$-matrix (see, for example \cite[Ch. IV, \S2]{CL}). But any formal
Laurent series $\tilde y=\sum_{j\geqslant-N}c_jz^j$, $c_j\in{\mathbb C}^n$, satisfying the system (\ref{System})
is in fact a Laurent polynomial, since relations
$$
A\,c_j=j\,c_j, \qquad j\geqslant-N,
$$
can hold only for a finite number of non-zero vectors $c_j$.

In the case when the singular point $z=0$ is {\it irregular} (i.~e., is not a regular singular one),
J.-P.~Ramis \cite{Ram} has suggested the following method for estimating the Gevrey order of a formal
power series solution. To the linear differential operator
$$
L=a_n(z)\delta^n+a_{n-1}(z)\delta^{n-1}+\ldots+a_0(z),
$$
which corresponds to the equation (\ref{LDEdelta}), one attaches its {\it Newton polygon} ${\cal N}(L)$, the
boundary curve of the smallest convex set containing the union of the sets
$$
X_i=\{(x,y)\in{\mathbb R}^2\,|\,x\leqslant i,\, y\geqslant{\rm ord}_0\,a_i(z)\},
\qquad i=0,1,\ldots,n.
$$
Let us note that in the Fuchsian case the Newton polygon consists of two line segments, a horizontal one and a
vertical one (see Pic. 1). But if $z=0$ is an irregular singular point of the equation (\ref{LDEdelta}), then the Newton polygon
contains line segments with positive slopes. Let $0<r_1<\ldots<r_m<\infty$ be all of such slopes ($m\leqslant n$, see Pic. 2).
Then, as Ramis's theorem asserts, {\it any formal power series satisfying the equation $(\ref{LDEdelta})$ has
the Gevrey type of precise order $s\in\{0,1/r_1,\ldots,1/r_m\}$}.

As a generalization of Ramis's theorem one regards the following theorem for a non-linear differential equation
\begin{eqnarray}\label{ODE}
F(z,u,\delta u,\ldots,\delta^nu)=0,
\end{eqnarray}
where $F(z,y_0,y_1,\ldots,y_n)$ is a holomorphic function in a neighbourhood of $0\in{\mathbb C}^{n+2}$.


\medskip
\begin{figure}[h!] \begin{picture}(240, 120)
\put(10, 1){\includegraphics*[height=4cm]{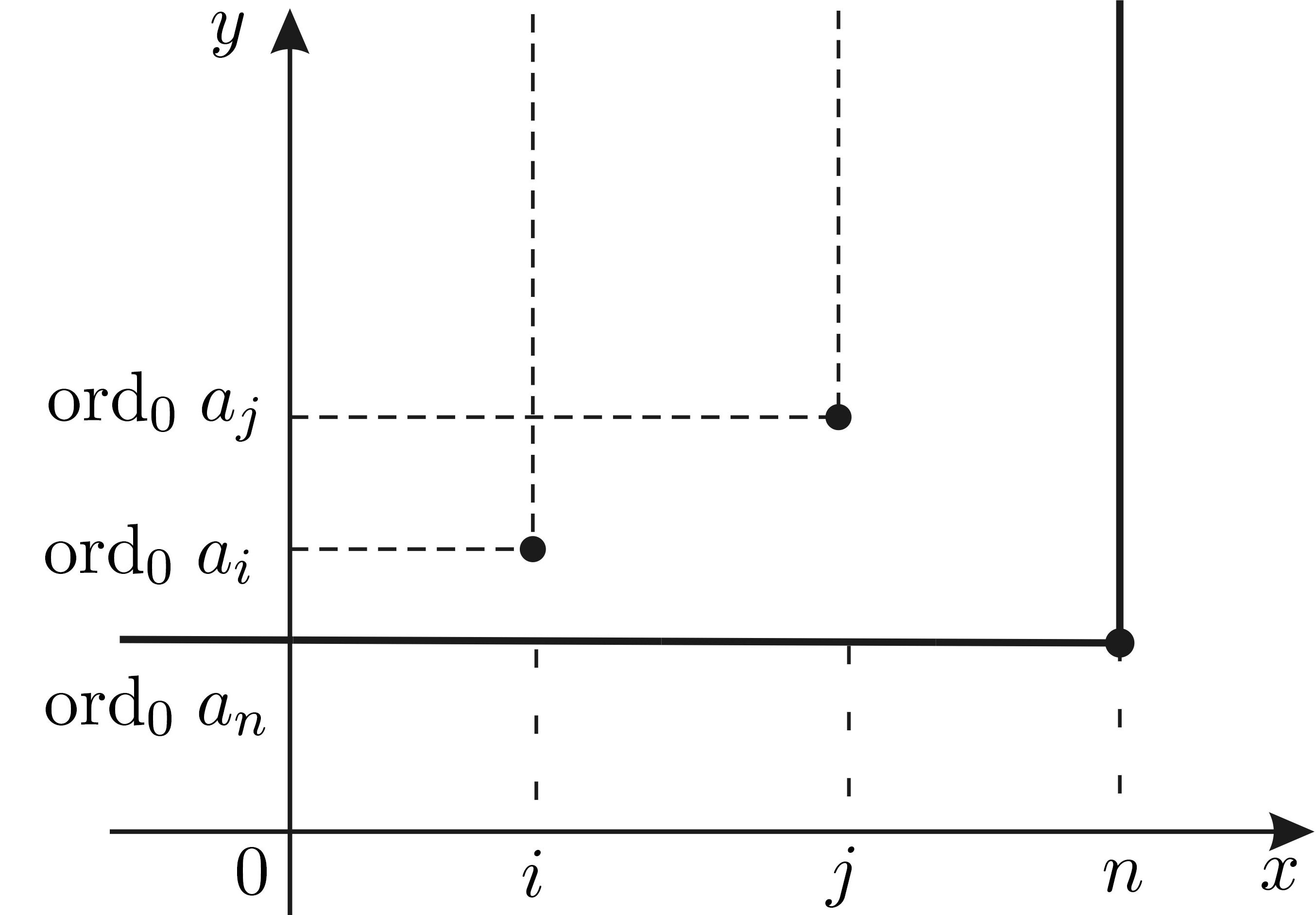}}
\put(270, 0){\includegraphics*[height=4cm]{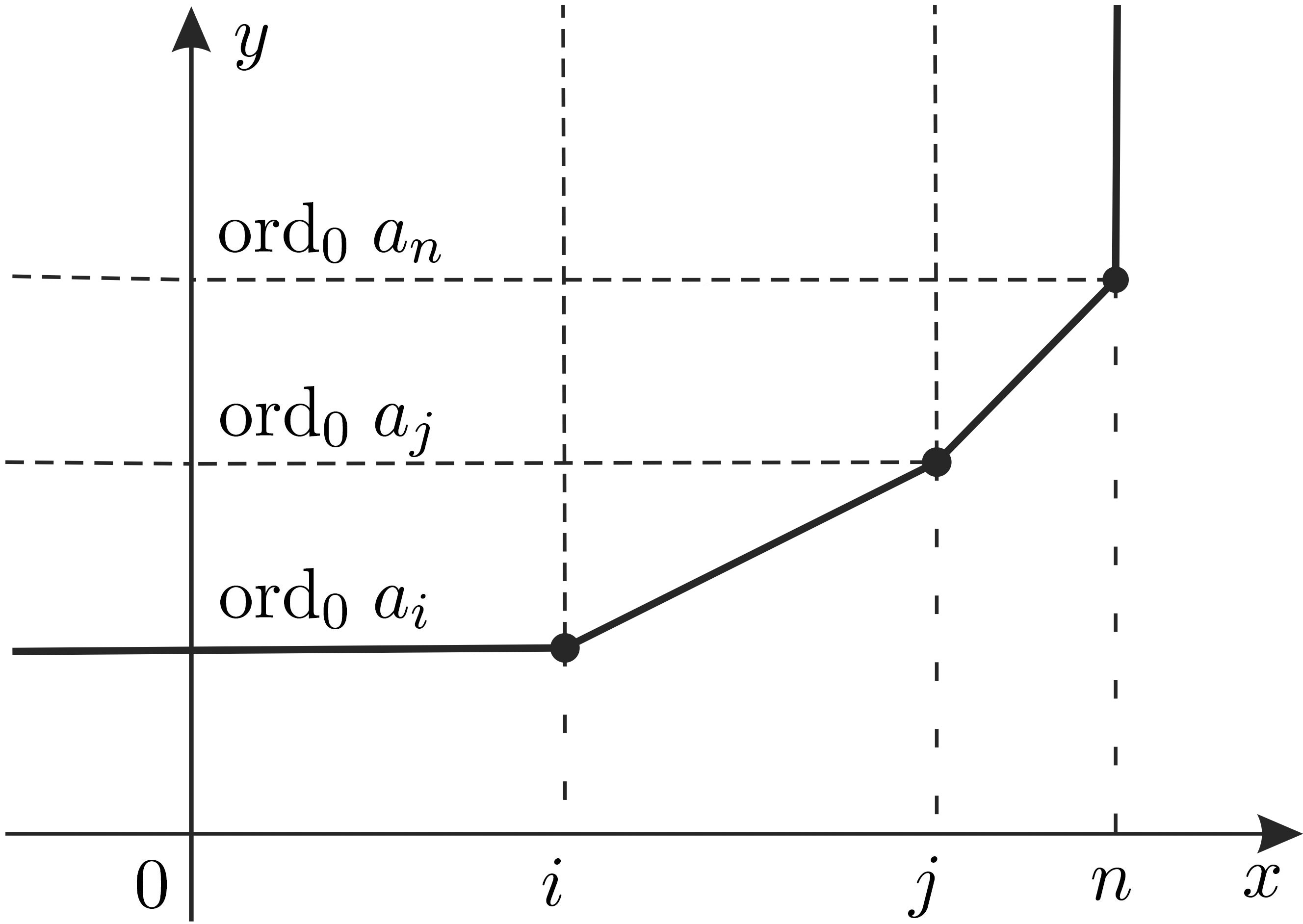}}
\end{picture}
\begin{center}
\begin{tabular}{lll}
\small Pic. 1.\quad The Newton polygon when $z=0$ is & \qquad\phantom{.} & \small Pic. 2.\quad The Newton polygon when $z=0$ is \\
\small a Fuchsian singular point of the operator $L$. &\qquad\phantom{.} & \small an irregular singular point of the operator $L$.\\
\quad& \qquad \phantom{.}&      \small  Here $\displaystyle r_1=({\rm ord}_0\,a_j-{\rm ord}_0\,a_i)/(j-i),\;$  \\
\quad&\qquad\phantom{.}&\small $\displaystyle r_2=({\rm ord}_0\,a_n-{\rm ord}_0\,a_j)/(n-j).$
\end{tabular}
\end{center}
\end{figure}

\medskip

{\bf Theorem 1} (B.\,Malgrange \cite{Mal}). {\it Let a formal series $\hat\varphi\in{\mathbb C}[[z]]$,
$\hat\varphi(0)=0$, satisfy the equation $(\ref{ODE})$ $($that is, $F(z, \Phi)=F(z,\hat\varphi,\delta\hat\varphi,\ldots,
\delta^n\hat\varphi)=0)\,$ and $\,\displaystyle\frac{\partial F}{\partial y_n}(z,\Phi)\ne0$. Then}

a) {\it if $z=0$ is a Fuchsian singular point of the operator $\displaystyle L_{\hat\varphi}=
\sum_{i=0}^n\frac{\partial F}{\partial y_i}(z,\Phi)\delta^i$, then the series $\hat\varphi$
converges in a neighbourhood of zero;}

b) {\it if $z=0$ is an irregular singular point of the operator $L_{\hat\varphi}$, and $r$ is the smallest of
positive slopes of the Newton polygon ${\cal N}(L_{\hat\varphi})$, then the formal power series $\hat\varphi$
has the Gevrey type of order $s=1/r$.}
\medskip

The part b) of Malgrange's theorem has been precised by Y.~Sibuya \cite[App. 2]{Si} as follows: {\it the formal power
series $\hat\varphi$ has the Gevrey type of the precise order $s\in\{0,1/r_1,\ldots,1/r_m\}$, where
$0<r_1<\ldots<r_m<\infty$ are all of positive slopes of the Newton polygon ${\cal N}(L_{\hat\varphi})$.}

The Malgrange--Sibuya theorem (Theorem 1a) on a sufficient condition of the convergence of a formal solution of an ordinary differential
equation has been proved by these authors in the different ways. First the equation (\ref{ODE}) is transformed to some special form,
then Malgrange uses the theorem on an implicit mapping for Banach spaces, while Sibuya applies the fundamental Ramis--Sibuya theorem \cite{RS}
on asymptotic expansions. We expose Malgrange's proof in details in the next section, then in the Section 3 we give an analytic proof which
allows to estimate the radius of convergence of the power series $\hat\varphi$ (Malgrange's and Sibuya's theorems do not contain estimates
for the radius of convergence). An idea of our proof is based on the construction of a majorant equation and was already appeared in the
article \cite[Ch. 1, \S7]{BG}.
\\

{\bf\S2. Malgrange's proof}
\\

For each natural $k$ the formal power series $\hat\varphi$ may be represented in the form
$$
\hat\varphi=\varphi_k+z^k\hat\psi, \qquad \hat\psi(0)=0.
$$

{\bf Lemma 1.} {\it For a sufficiently large $k$ $($under the assumptions of Theorem 1a$)$, the formal power series $\hat\psi$
satisfies the relation
$$
\overline{L}(\delta+k)\hat\psi=zM(z,\hat\psi,\delta\hat\psi,\ldots,\delta^n\hat\psi),
$$
where $\overline{L}$ is a polynomial of degree $n$, $M$ is a holomorphic function in a neighbourhood of $0\in{\mathbb C}^{n+2}$.}
\medskip

{\bf Proof.} Taking into consideration the equality $\delta(z^k\hat\psi)=z^k(\delta+k)\hat\psi$ we have the relations
$$
\delta^i(z^k\hat\psi)=z^k(\delta+k)^i\hat\psi, \qquad i=1,\ldots,n,
$$
therefore,
$$
\Phi=(\varphi_k,\delta\varphi_k,\ldots,\delta^n\varphi_k)+
z^k(\hat\psi,(\delta+k)\hat\psi,\ldots,(\delta+k)^n\hat\psi)=\Phi_k+z^k\Psi.
$$
Further applying Taylor's formula we obtain
\begin{eqnarray}\label{Taylor}
0=F(z,\Phi_k+z^k\Psi)&=&F(z,\Phi_k)+z^k\sum_{i=0}^n\frac{\partial F}{\partial y_i}(z,\Phi_k)
(\delta+k)^i\hat\psi+ \nonumber \\ & &+z^{2k}\sum_{p,q=0}^nH_{pq}(z,\Phi_k,z^k\Psi)
(\delta+k)^p\hat\psi(\delta+k)^q\hat\psi,
\end{eqnarray}
where $H_{pq}$ are holomorphic functions in a neighbourhood of $0\in{\mathbb C}^{2n+3}$.

Let $\displaystyle l={\rm ord}_0\frac{\partial F}{\partial y_n}(z,\Phi)$. Then by the assumption of Theorem 1a,
$$
{\rm ord}_0\frac{\partial F}{\partial y_i}(z,\Phi)\geqslant l, \qquad i=0,1,\ldots,n.
$$
Let $b_iz^l$ be a summand that the formal power series $\displaystyle\frac{\partial F}{\partial y_i}(z,\Phi)$ begins with:
$$
\frac{\partial F}{\partial y_i}(z,\Phi)=b_iz^l+\ldots,
\quad i=0,1,\dots,n \quad (b_i\in\mathbb C, \;b_n\neq 0).
$$
Define a polynomial
\begin{equation}
\overline{L}(\xi)=\sum_{i=0}^nb_i\xi^i\label{Euler}
\end{equation}
of degree $n$ and choose a number $k_0$ such that for any natural $k>k_0$ the inequality $\overline{L}(k)\ne0$ holds.
Now we show that the number $k$ from the statement of the lemma can be taken as $k=\max(k_0,l+1)$.

Let us note that
$$
{\rm ord}_0\left(\frac{\partial F}{\partial y_i}(z,\Phi)-
\frac{\partial F}{\partial y_i}(z,\Phi_k)\right)\geqslant l+1, \qquad i=0,1,\ldots,n
$$
(this follows from Taylor's formula applied to this difference), therefore,
$\displaystyle \frac{\partial F}{\partial y_i}(z,\Phi_k)=b_iz^l+o(z^l)$ in a neighbourhood of zero. From the relations 
(\ref{Taylor}) and $\hat\psi(0)=0$ it follows that
$$
{\rm ord}_0F(z,\Phi_k)\geqslant k+l+1.
$$
Hence the relation (\ref{Taylor}) can be divided by $z^{k+l}$, and we obtain the equality of the form
$$
\overline{L}(\delta+k)\hat\psi-zM(z,\hat\psi,\delta\hat\psi,\ldots,\delta^n\hat\psi)=0,
$$
where the polynomial $\overline{L}$ is defined by the formula (\ref{Euler}), $M$ is a holomorphic function in a neighbourhood of
$0\in{\mathbb C}^{n+2}$.  $\quad\Box$
\medskip

Thus the formal power series $\hat\psi=\sum_{j=1}^\infty c_jz^j$ is a solution of the ordinary differential equation
\begin{eqnarray}\label{ODE2}
\overline{L}(\delta+k)v=zM(z,v,\delta v,\ldots,\delta^nv).
\end{eqnarray}
Let us prove the convergence of this series in some neighbourhood of zero using the theorem on an implicit mapping for Banach spaces. 
We recall a part of this theorem that we need here (see, for example, \cite[Th. 10.2.1]{Die}).
\medskip

{\it Let $E$, $F$, $G$ be Banach spaces, $A$ an open subset of the direct product $E\times F$,
$f:\,A\longrightarrow G$ a continuously differentiable mapping. Consider a point $(x_0,y_0)\in A$ such that
$f(x_0,y_0)=0$ and $\displaystyle \frac{\partial f}{\partial y}(x_0,y_0)$ is a} bijective {\it linear mapping from $F$ to $G$.

Then there is a neighbourhood $U_0\subset E$ of the point $x_0$ and a unique continuous mapping $u:\,U_0\longrightarrow F$
such that $u(x_0)=y_0$, $(x,u(x))\in A$ and $f(x, u(x))=0$ for any $x\in U_0$.}
\medskip

Note that in view of the condition $\overline L(j+k)\ne0$, $j=1,2,\ldots$, the coefficient $c_1$ and all the other coefficients
of the power series $\hat\psi$ are uniquely determined (every coefficient is expressed via the previous ones). Thus, one has a unique 
power series that satisfies the equation (\ref{ODE2}).

For each $m=0,1,\ldots,n$ let us define a Banach space
$$
\textstyle H^m=\Bigl\{\psi=\sum\limits_{j\geqslant1} a_jz^j\;\;\Bigl|\Bigr.\; \sum\limits_{j\geqslant1} j^m|a_j|<\infty\Bigr\}
$$
with the corresponding norm $\|\psi\|_m=\sum_{j\geqslant1} j^m|a_j|=\|\delta^m\psi\|_0$. It is not difficult to see that 
$H^m\subset H^{m-1}$, $m=1,2,\ldots,n$, and all the spaces $H^m$ are contained in the space of functions holomorphic in the open 
unit disk $D_1$ and contain the space of functions holomorphic in the closed unit disk $\overline D_1$:
$$
{\cal O}(\overline D_1)\subset H^n\subset\ldots\subset H^1\subset H^0\subset{\cal O}(D_1).
$$

Let the function $M(z,y_0,y_1,\ldots,y_n)$ in the equation (\ref{ODE2}) be holomorphic in a polydisk $\{|z|<\varepsilon,
|y_0|<\varepsilon,\ldots,|y_n|<\varepsilon\}$. In the product ${\mathbb C}\times H^n$ of the Banach spaces we consider an open subset
$$
A=\{(\lambda,\psi):\,|\lambda|<\varepsilon, \|\psi\|_n<\varepsilon\}
$$
and define a continuously differentiable mapping
$$
f: (\lambda,\psi)\mapsto\overline{L}(\delta+k)\psi-
\lambda zM(\lambda z,\psi,\delta\psi,\ldots,\delta^n\psi)
$$
from $A$ to $H^0$, with $f(0,0)=0$. Let us show that the linear operator
$$
\frac{\partial f}{\partial\psi}(0,0)=\overline{L}(\delta+k): H^n\longrightarrow H^0
$$
is bijective. Indeed,
$$
\overline{L}(\delta+k)\,a_jz^j=a_j\overline{L}(j+k)z^j=0\;\Longleftrightarrow\;a_j=0,
$$
therefore, ${\rm ker}\,\overline{L}(\delta+k)=\{0\}$. In the same time, if $\sum_{j=1}^{\infty}a_jz^j\in H^0$,
then $\sum_{j=1}^{\infty}(a_j/\overline{L}(j+k))z^j\in H^n$, that is, the image of the operator $\overline{L}(\delta+k)$
coincides with $H^0$.

Thus, the mapping $f$ satisfies the theorem on an implicit mapping, hence there are a real number $\nu>0$ and a function  
$\psi_{\nu}\in H^n\subset{\cal O}(D_1)$ such that
$$
\overline{L}(\delta+k)\psi_{\nu}(z)-
\nu zM(\nu z,\psi_{\nu}(z),\delta\psi_{\nu}(z),\ldots,\delta^n\psi_{\nu}(z))=0.
$$
But then the function $\psi_{\nu}(z/\nu)\in{\cal O}(D_{\nu})$ is a solution of the equation (\ref{ODE2}), and the power series 
$\hat\psi$ converges in the disk $D_{\nu}$ of radius $\nu$. (Here we use the relations
$(\delta^p\psi_{\nu})(z/\nu)=\delta^p(\psi_{\nu}(z/\nu))$, $p=1,2,\ldots,n$.)
\medskip

Let us note that in the proof above it is essentially that $\overline L$ is a polynomial of degree $n$. If the degree would be 
less than $n$ (this is the case, when $z=0$ is an irregular singular point of the operator $L_{\hat \varphi}$), for example, would be
equal to $n-1$, then the linear operator $\overline{L}(\delta+k): H^n\longrightarrow H^0$
is not surjective, since its image coincides with $H^1$, and we cannot apply the theorem on an implicit mapping.
\\
\\

{\bf\S3. Proof by the majorant method. An estimate for the radius of convergence}
\\

Here we construct a differential equation majorant for the equation (\ref{ODE2}), in the sense that
it will have a unique solution $\psi=\sum_{j=1}^{\infty}C_jz^j$ ($C_j\geqslant0$) holomorphic in a neighbourhood of zero, 
and its power series will be majorant for the formal power series $\hat\psi$. Besides that, we will obtain an estimate for 
the radius of convergence of the power series $\psi$ and, hence, power series $\hat\psi$.

Let us write the function $M(z,y_0,y_1,\ldots,y_n)$ from the right hand side of the equation (\ref{ODE2}) as a power series
converging in some neighbourhood of $0\in{\mathbb C}^{n+2}$:
$$
M(z,y_0,y_1,\ldots,y_n)=\sum_{p=0}^{\infty}\sum_{{\bf q}=(q_0,q_1,\ldots,q_n)\in{\mathbb Z}_+^{n+1}}
\alpha_{p,{\bf q}}\,z^p\,y_0^{q_0}y_1^{q_1}\ldots y_n^{q_n}, \quad \alpha_{p,{\bf q}}\in\mathbb C.
$$
As we noted before, each coefficient $c_j$ of the formal power series $\hat\psi=\sum_{j=1}^{\infty}c_jz^j$ is uniquely determined 
by the previous ones. Now we find exact expressions for these coefficients using the relation
$$
\overline{L}(\delta+k)\hat\psi=zM(z,\hat\psi,\delta\hat\psi,\ldots,\delta^n\hat\psi).
$$

Denoting by $\widehat F(z)$ the formal power series $M(z,\hat\psi,\delta\hat\psi,\ldots,\delta^n\hat\psi)$ one has
$$
\overline{L}(j+k)c_j=\frac{\widehat F^{(j-1)}(0)}{(j-1)!}, \qquad j=1,2,\ldots\;.
$$
In order to express $\widehat F^{(j-1)}(0)$ we use a formula for the derivative of a product:
$$
(f_1\ldots f_q)^{(m)}=\sum_{m_1+\ldots+m_q=m}\frac{m!}{m_1!\ldots m_q!}f_1^{(m_1)}\ldots f_q^{(m_q)},
\qquad f_i\in{\mathbb C}[[z]], \quad m=1,2,\ldots\;.
$$
For $y_i=\delta^i\hat\psi=\sum_{j=1}^{\infty}j^i\,c_jz^j$ ($i=0,1,\ldots,n$) we obtain for each $q\leqslant m$:
$$
(y_i^q)^{(m)}(0)=\sum_{m_1+\ldots+m_q=m}m!\,\frac{y_i^{(m_1)}(0)}{m_1!}\ldots\frac{y_i^{(m_q)}(0)}{m_q!}=
\sum_{m_1+\ldots+m_q=m}m!(m_1\ldots m_q)^i\,c_{m_1}\ldots c_{m_q}
$$
($=0$, if $q>m$), and also
$$
\Bigl(z^p\,y_0^{q_0}y_1^{q_1}\ldots y_n^{q_n}\Bigr)^{(j-1)}(0)=\sum_{p+j_0+\ldots+j_n=j-1}\frac{(j-1)!}{j_0!j_1!\ldots j_n!}
(y_0^{q_0})^{(j_0)}(y_1^{q_1})^{(j_1)}\ldots(y_n^{q_n})^{(j_n)}\Bigl|_{z=0}.
$$
Therefore,
\begin{eqnarray}\label{f1}
\frac{\widehat F^{(j-1)}(0)}{(j-1)!}=\sum_{p+q_0+\ldots+q_n\leqslant j-1}\alpha_{p,{\bf q}}\sum_{j_0+\ldots+j_n=j-1-p}
\frac{(y_0^{q_0})^{(j_0)}}{j_0!}\frac{(y_1^{q_1})^{(j_1)}}{j_1!}\ldots\frac{(y_n^{q_n})^{(j_n)}}{j_n!}\Bigl|_{z=0},
\end{eqnarray}
where
\begin{eqnarray}\label{f2}
\frac{(y_i^{q_i})^{(j_i)}(0)}{j_i!}=\sum_{m_1+\ldots+m_{q_i}=j_i}(m_1\ldots m_{q_i})^i\,c_{m_1}\ldots c_{m_{q_i}},
\qquad i=0,1,\ldots,n.
\end{eqnarray}
Thus,
\begin{eqnarray}\label{c_rel}
\overline{L}(1+k)c_1=\alpha_{0,0}, \qquad\overline{L}(j+k)c_j=P_j(c_1,\ldots,c_{j-1},\{\alpha_{p,{\bf q}}\}),
\quad j=2,3,\ldots\;,
\end{eqnarray}
where $P_j$ is a polynomial with positive coefficients that is restored by the formulas (\ref{f1}), (\ref{f2}).

Now we consider an equation majorant, as will be shown further, for the differential equation (\ref{ODE2}):
\begin{eqnarray}\label{ODE3}
\sigma\,\delta^nv=z\widetilde M (z,\delta^nv), \qquad\sigma>0,
\end{eqnarray}
where
$$
\widetilde M (z,w)=\sum_{p=0}^{\infty}\sum_{{\bf q}=(q_0,q_1,\ldots,q_n)\in{\mathbb Z}_+^{n+1}}
|\alpha_{p,{\bf q}}|\,z^p\,w^{q_0}w^{q_1}\ldots w^{q_n}
$$
is a holomorphic function in a neighbourhood of the point $(0,0)\in{\mathbb C}^2$, and its power series expansion 
is produced from the power series expansion of the function $M(z,y_0,y_1,\ldots,y_n)$ by changing all the coefficients 
$\alpha_{p,{\bf q}}$ to their absolute values and the variables $y_0,y_1,\ldots,y_n$ to one variable $w$. The value 
$\sigma$ is defined by the formula
$$
\sigma=\inf_{j\in\mathbb N}\frac{|\overline L(j+k)|}{j^n}.
$$
This value is positive, as $\overline L(j+k)\ne0$ for $j\in\mathbb N$ and
$\lim\limits_{j\rightarrow\infty}|\overline L(j+k)|/j^n=|b_n|>0$ (recall that $\overline L(\xi)=\sum_{i=0}^nb_i\xi^i$).
\medskip

{\bf Lemma 2.} {\it The differential equation $(\ref{ODE3})$ has a unique solution of the form
$\psi=\sum_{j=1}^{\infty}C_jz^j$ $($i.~e., $\psi(0)=0)$ holomorphic in a neighbourhood of zero, and its power series is 
majorant for the power series $\hat\psi=\sum_{j=1}^{\infty}c_jz^j$ satisfying the equation $(\ref{ODE2})$.}
\medskip

{\bf Proof.} Existence and uniqueness of a solution $\psi=\sum_{j=1}^{\infty}C_jz^j$ of the differential equation $(\ref{ODE3})$ 
follows from existence and uniqueness of a holomorphic in a neighbourhood of zero solution $w=w(z)=\sum_{j=1}^{\infty}a_jz^j$  
($w(0)=0$) of the equation $\sigma\,w=z\widetilde M(z,w)$, as the latter satisfies the assumptions of Cauchy's theorem on an implicit 
function at the point  $(0,0)\in{\mathbb C}^2$ (then $C_j=a_j/j^n$).

Let us show that the power series $\psi=\sum_{j=1}^{\infty}C_jz^j$ is majorant for the power series $\hat\psi=\sum_{j=1}^{\infty}c_jz^j$. 
We will find the coefficients $C_j$ using the relation
$$
\sigma\,\delta^n\psi=z\widetilde M(z,\delta^n\psi).
$$
Denoting by $\widetilde F(z)$ the function $\widetilde M(z,\delta^n\psi)$ holomorphic in a neighbourhood of zero, one has
$$
\sigma\,j^nC_j=\frac{\widetilde F^{(j-1)}(0)}{(j-1)!}, \qquad j=1,2,\ldots\;.
$$
Keeping in mind the difference between $\widetilde F(z)$ and $\widehat F(z)$, similarly to the formulas (\ref{f1}), (\ref{f2})
for $\widehat F^{(j-1)}(0)/(j-1)!$ we obtain
\begin{eqnarray}\label{f3}
\frac{\widetilde F^{(j-1)}(0)}{(j-1)!}=\sum_{p+q_0+\ldots+q_n\leqslant j-1}|\alpha_{p,{\bf q}}|\sum_{j_0+\ldots+j_n=j-1-p}
\frac{(w^{q_0})^{(j_0)}}{j_0!}\frac{(w^{q_1})^{(j_1)}}{j_1!}\ldots\frac{(w^{q_n})^{(j_n)}}{j_n!}\Bigl|_{z=0},
\end{eqnarray}
where
\begin{eqnarray}\label{f4}
\frac{(w^{q_i})^{(j_i)}(0)}{j_i!}=\sum_{m_1+\ldots+m_{q_i}=j_i}(m_1\ldots m_{q_i})^n\,C_{m_1}\ldots C_{m_{q_i}},
\qquad i=0,1,\ldots,n.
\end{eqnarray}
Hence,
\begin{eqnarray}\label{C_rel}
\sigma\,C_1=|\alpha_{0,0}|\in{\mathbb R}_+, \qquad\sigma\,j^nC_j=
\widetilde P_j(C_1,\ldots,C_{j-1},\{|\alpha_{p,{\bf q}}|\})\in{\mathbb R}_+,\quad j=2,3,\ldots\;,
\end{eqnarray}
where $\widetilde P_j$ is a polynomial with positive coefficients that is restored by the formulas (\ref{f3}), (\ref{f4}). 
Let us note that by the construction one has
$$
\widetilde P_j(C_1,\ldots,C_{j-1},\{|\alpha_{p,{\bf q}}|\})\geqslant
P_j(C_1,\ldots,C_{j-1},\{|\alpha_{p,{\bf q}}|\})
$$
(compare the formulas (\ref{f2}) and (\ref{f4})).

Using the obtained recurrence relations (\ref{c_rel}) and (\ref{C_rel}) for the coefficients
$c_j$ and $C_j$ we see that
$$
|\overline L(1+k)||c_1|=|\alpha_{0,0}|=\sigma\,C_1\quad\Longrightarrow\quad
|c_1|=\frac{\sigma}{|\overline L(1+k)|}\,C_1\leqslant C_1,
$$
and finish the proof by the induction (the second inequality below):
\begin{eqnarray*}
|\overline L(j+k)||c_j| & = & |P_j(c_1,\ldots,c_{j-1},\{\alpha_{p,{\bf q}}\})|\leqslant
P_j(|c_1|,\ldots,|c_{j-1}|,\{|\alpha_{p,{\bf q}}|\})\leqslant\\
& \leqslant & P_j(C_1,\ldots,C_{j-1},\{|\alpha_{p,{\bf q}}|\})\leqslant
\widetilde P_j(C_1,\ldots,C_{j-1},\{|\alpha_{p,{\bf q}}|\})=\\
& = & \sigma\,j^nC_j\quad\Longrightarrow\quad
|c_j|\leqslant\frac{\sigma\,j^n}{|\overline L(j+k)|}\,C_j\leqslant C_j, \quad j=2,3,\ldots\;.
\end{eqnarray*} 
The proof is finished.
\medskip

{\bf Proposition 1.} {\it Let the function $M(z,y_0,y_1,\ldots,y_n)$ from the right hand side of the equation $(\ref{ODE2})$
be holomorphic in a neighbourhood of a closed  polydisk 
$$
\overline\Delta=\{|z|\leqslant r, |y_0|\leqslant\rho,\ldots,
|y_n|\leqslant\rho\},\;\; \mu=\max\limits_{\overline\Delta}|M|.
$$ 
Then the power series $\hat\psi=\sum_{j=1}^{\infty}c_jz^j$ satisfying the equation $(\ref{ODE2})$ converges in the disk}
$$
D_R=\Bigl\{|z|<r\,\frac{\rho}{\rho+\mu r/\sigma N}\Bigr\}, \qquad N=(n+1)^{n+1}/(n+2)^{n+2}.
$$

{\bf Proof.} It is enough to check that the power series $\psi=\sum_{j=1}^{\infty}C_jz^j$, satisfying the equation (\ref{ODE3})
and majorant for $\hat\psi$, converges in the disk $D_R$. It will follow from the convergence of the power series 
$\delta^n\psi=\sum_{j=1}^{\infty}j^nC_jz^j$, which is a solution of the equation $\sigma\,w=z\widetilde M(z,w)$. Here we apply 
the majorant method from the proof of the theorem on an implicit function (see, for example, \cite[Ch. IX, \S193]{Gou}).

As
$$
\widetilde M (z,w)=\sum_{p=0}^{\infty}\sum_{{\bf q}=(q_0,q_1,\ldots,q_n)\in{\mathbb Z}_+^{n+1}}
|\alpha_{p,{\bf q}}|\,z^p\,w^{q_0}w^{q_1}\ldots w^{q_n},\quad
|\alpha_{p,{\bf q}}|\leqslant\frac{\mu}{r^p\rho^{q_0}\rho^{q_1}\ldots\rho^{q_n}},
$$
the power series $\delta^n\psi$ satisfying the equation $\sigma\,w=z\widetilde M(z,w)$ converges apriori in a disk where the solution 
(equal to zero for $z=0$) of the equation 
$$
\sigma\,w=z\mu\sum_{p=0}^{\infty}\sum_{{\bf q}=(q_0,q_1,\ldots,q_n)\in{\mathbb Z}_+^{n+1}}
\Bigl(\frac zr\Bigr)^p\Bigl(\frac w{\rho}\Bigr)^{q_0}\Bigl(\frac w{\rho}\Bigr)^{q_1}\ldots
\Bigl(\frac w{\rho}\Bigr)^{q_n}=\frac{\mu z}{(1-z/r)(1-w/\rho)^{n+1}}
$$
is holomorphic. Let us rewrite the last equation in the form
$$
f(z,w)=\sigma\,w(1-w/\rho)^{n+1}-\frac{\mu z}{1-z/r}=0.
$$
By the theorem on an implicit function, its solution admits a single valued branch that is equal to zero for $z=0$ and holomorphic in 
a disk $\{|z|<|z_0|\}$, where $z_0$ is determined from the system of the equations $f(z,w)=0$,
$\frac{\partial f(z,w)}{\partial w}=0$. A unique solution $w_0=\rho/(n+2)$, $z_0=r\,\frac{\rho}{\rho+\mu r/\sigma N}$ of this system leads 
to an estimate for the radius of convergence of the power series $\delta^n\psi$, hence of $\psi$ and $\hat\psi$ as well. $\quad \Box$

\bigskip
\bigskip

R.\,R.\,Gontsov, The Institute for Information Transmission Problems of RAS,

Moscow, Russia, gontsovrr@gmail.com
\\

I.\,V.\,Goryuchkina, The Keldysh Institute of Applied Mathematics of RAS,

Moscow, Russia, igoryuchkina@gmail.com 

\end{document}